\documentclass[11pt]{amsart}
\usepackage{epsfig,amssymb,latexsym,amscd,xypic}

\theoremstyle{plain}

\newtheorem{theo}{Theorem}
\newtheorem{coro}{Corollary}
\newtheorem{lema}{Lemma}

\theoremstyle{remark}

\newtheorem{example}{Example}

\theoremstyle{definition}

\newtheorem{defi}{Definition}
\newtheorem{obse}{Observation}

\newcommand{\qqed}
                {
                \qed\vspace*{3mm}
                }

\newcommand{\qedennonce}
                {
                \ \vspace*{-10pt} \\
                \qed
                \vspace*{3mm}
                }

\newcommand{\dem}
                {
                {\sl {\sc Proof.}}
                }

\newcommand{\codim}
                {
                \operatorname{codim}
                }

\begin{document}
\title{Generalized Cayley's $\Omega$-processes} 
\author{Walter Ferrer Santos, Alvaro Rittatore}
\date{\today}
\thanks{The first author would like to thank, Csic-UDELAR and
  Conicyt-MEC, the second author would like to thank FCE-MEC}

\begin{abstract}
In this paper we generalize some constructions and results 
due to Cayley and Hilbert. We define the concept of $\Omega$--process
for an arbitrary algebraic monoid with zero and unit group $G$. Then  we
show how to produce from the process and 
for a linear rational representation of $G$, 
a number of elements
of the ring of $G$-invariants, that is large
enough as to guarantee its finite generation. Moreover, we give an explicit
construction of all $\Omega$-processes for general reductive monoids
and, in the case of the monoid of all the $n^2$ matrices, compare our
construction with Cayley's definition.
\end{abstract}

\maketitle

\section{Introduction}

We assume throughout that the base field  $\Bbbk$ is 
algebraically closed of characteristic zero and that 
all the geometric and
algebraic objetcs are defined over $\Bbbk$. 
A {\em linear algebraic monoid}\/ is an affine normal
algebraic variety $M$ with an associative  product $M\times M\to M$ which
is a morphism of algebraic $\Bbbk$-varieties, and that is also equipped with 
a neutral element $1\in M$ for this
product. In this situation, it can be proved that the group of invertible
elements of $M$ -- called the unit group
of $M$ --, $G(M)=\{
g\in M \mathrel{|} \,\,\,\exists \,\, g^{-1}\}$ is an affine
algebraic group, that is open in $M$.  
A {\em reductive monoid}\/ is an irreducible normal
linear algebraic  monoid  whose unit group is reductive.

In this work we define an $\Omega$-process for a general affine
monoid with $0$ generalizing the classical definition of Cayley for
$\operatorname {M}_n(\Bbbk)$. 
An $\Omega$-process for a monoid $M$ is  a  left and right semi-equivariant
map $\Omega:\Bbbk[M] \rightarrow 
\Bbbk[M]$ (see Section \ref{section:process} for the precise
definition).  We also show 
how to construct from $\Omega$ a
family of Reynolds operators on the polynomial representations of $M$ and
more importantly -- following Hilbert's methods as presented in
\cite{kn:Hil} -- we  use $\Omega$ to prove the finite generation
of the invariants of the linear representations of $G(M)$.
Finally, given a reductive monoid $M$, we characterize the family of all
$\Omega$-processes for $M$. 

For the basic notations of the theory of affine algebraic groups and
in particular for the notion of Reynolds operators, 
integrals and invariants see \cite{kn:nosotros}. For a review of know
results about reductive monoids and their representations, see
\cite{kn:rennersurvey}.

\section{Preliminaries}

In a similar way than for the case of affine algebraic groups, if $M$
is an affine algebraic monoid, the multiplication of $M$ induces a
morphism of algebras $\Delta: \Bbbk [M] \rightarrow \Bbbk [M] \otimes
\Bbbk [M]$. This map, together with the evaluation at the neutral
element $\varepsilon: \Bbbk [M] \rightarrow \Bbbk$ induces a bialgebra
structure on $\Bbbk [M]$. Moreover, if the monoid $M$ is endowed with a zero
element, i.e.\ an element $0 \in M$ such that $0m=m0=0$ for all $m \in
M$, 
then the evaluation at $0$ induces an algebra morphism $\nu: \Bbbk [M]
\rightarrow \Bbbk$ such that $u\nu \mathrel{\star} \operatorname {Id} = u\nu =
\operatorname {Id} \mathrel{\star} u\nu$, where $u:\Bbbk \rightarrow \Bbbk [M]$
is the unit of the algebra $\Bbbk [M]$ and $\star$ denotes the
convolution product on the endomorphisms of $\Bbbk[M]$: if $T,S \in
\operatorname {End}\bigl(\Bbbk[M]\bigr)$, 
then $T \mathrel{\star} S \in \operatorname {End}\bigl(\Bbbk[M]\bigr)$ is
defined as $(T\star 
S)(f)=\sum T(f_1)S(f_2)$ if $\Delta(f)=\sum f_1 \otimes f_2$.
  
In particular if $M$ and $N$ are affine algebraic monoids, a
homomorphism of algebraic monoids $\Phi: M \rightarrow N$, 
is a morphism of algebraic
varieties that preserves de product and the unit element of the
monoids. Moreover, if the monoids have zero we assume that
$\Phi(0)=0$. One can easily show in this situation that the associated
map $\Phi^*:\Bbbk[N] \rightarrow \Bbbk[M]$ is a morphism of bialgebras.

The following result showing that an arbitrary affine monoid can be
linearized is easy to prove (see  \cite{kn:putcha-monoids} for a proof):

\begin{theo}
\label{theo:closedinm}
Let $M$ be an affine algebraic monoid and denote as
${\operatorname M}_n(\Bbbk)$ the monoid of $n\times n$-matrices with coefficients
in $\Bbbk$. Then there exist $n\in\mathbb
N$  and a closed inmersion $\varphi:M\to {\operatorname M}_n(\Bbbk)$, such that $\varphi\bigl(G(M)\bigr)=\varphi(M)\cap
\operatorname{GL}_n(\Bbbk)$. Moreover, if $M$ has a zero, then
$\varphi$ can be chosen in such a way that $\varphi(0)=0$. \qed
\end{theo}

Next, we recall some basic facts about reductive monoids
that will be used later.

Recall that any  algebraic monoid $M$ with unit group $G$ supports an
action of  $G\times G$, given by $(a,b)\cdot m=amb^{-1}$,
and this action has $G$ as an 
open orbit, see \cite{kn:am}. Thus,  reductive monoids are affine
spherical varieties and as such can be classified in combinatorial terms (see
for example \cite{kn:knop},\cite{kn:am}).

The following notations will be in force in this article.

Let $M$ be an algebraic monoid, and assume that $G=G(M)$ is a
reductive group.  We call $T$ a maximal torus of 
$G$, $B$ a Borel subgroup containing $T$, $B^-$ its
opposite Borel subgroup, denote as $\mathcal {X} (T)$ the set of weights 
and $\mathcal X_+(T)$ the semigroup of dominant weights with  respect to
$B$. We  call $W$ the Weyl group associated to $T$, $C=C(G)$
the Weyl chamber associated to $(B,T)$, $\alpha_1,\ldots
,\alpha_l$ and $\omega_1,\ldots, \omega_l$ the simple roots and
fundamental weights associated to $(B,T)$ respectively. Recall that
$\mathcal X(T)$ is partially ordered by the relation $\lambda\leq \mu$ if
and only if $\mu-\lambda\in \mathbb Z^+(\alpha_1,\ldots,
\alpha_l)$\,.

Finally, we call $\mathcal X_*(T)$ the set of one parameter subgroups (1-PS)
of $T$\, and identify $\mathcal X_*(T)$ with $\mathcal X(T)$ by means of a
$W$-invariant scalar product $\langle \cdot,\cdot\rangle$, in such a
way that $\langle \omega_i,\alpha_j^\vee\rangle=\delta_{ij}$, were
$\alpha_i^\vee$ is the coroot associated to $\alpha_i$.

If $\mathcal Q(G)=\mathcal X_*(T)\otimes_{\mathbb Z}\mathbb
Q=\mathcal X(T)\otimes_{\mathbb Z}\mathbb Q$, then the reductive
monoids with unit group $G$ are in one to one correspondence with   
with the pairs $(\mathcal C,\mathcal F)$ where $\mathcal 
F$ is a  subset of the simple roots and $\mathcal C\subset \mathcal
Q(G)$ is a rational poliedral cone generated by  $\mathcal F$ and a
finite number of elements of $-C(G)\subset \mathcal X(T)$, such that
the cone generated by $\mathcal C$ and all the simple roots is
strictly convex  (\cite[\S 4.2]{kn:am}). 

This result follows from 
the classification theory of spherical varieties in terms of
colored fans. Under the previous identification, we will denote
as $\mathcal C^\vee\subset \mathcal Q(G)$ the dual of the cone $\mathcal C$.

\section{Representations and character of algebraic monoids}

The correct category of representations of an affine algebraic monoid $M$
is the  category of the so called polynomial $M$-modules. 

If $V$ is a finite dimensional $\Bbbk$-space, we say that an action
$M \times V \rightarrow V$ is polynomial if the associated map $M
\rightarrow \operatorname {End}_{\Bbbk}(V)$ is a morphism of affine
algebraic monoids, in particular in the case that $M$ is endowed with 
$0 \in M$, $0 \cdot v = v$ for all 
$v \in V$.

The action map $M \times V \rightarrow V$ induces a morphism of 
$\Bbbk$-spaces 
$\chi_V: V \rightarrow V \otimes \Bbbk [M]$, that is a counital
comodule structure on $V$. 

The relationship between the structure $\chi_V$ and the action of $M$
is given as follows: $\chi_V(v)=\sum v_0 \otimes v_1$ if and only if
for all $m \in M$, $m \cdot v = \sum v_1(m)v_0$.
In this paper we consider -- in general --
the actions on the left side, but it is clear that all
the concepts we deal with have a formulation for a right action.

If $V$ is a polynomial $M$-module, we will denote as $V^M$ the subspace of
$M$-invariants and it follows easily that $V^M=\{v \in V:
\chi_{V}(v)=v \otimes 1\}$.

It is easy to prove, in the same way than for affine algebraic groups
and their representations, that the above correspondence is an
isomorphism between the category of finite dimensional
polynomial $M$-modules and the category of finite dimensional counital 
comodules for the bialgebra $\Bbbk [M]$.

If we drop the hypothesis on the dimension, we obtain the concept of
general polynomial $M$-module that is simply an arbitrary counital
comodule $V$ for the bialgebra $\Bbbk [M]$; the comodule structure
in $V$ will be denoted as $\chi_{V}: V \rightarrow V \otimes
\Bbbk[M]$. 

This corresponds to 
a vector space $V$ endowed with a locally finite polynomial 
action of $M$. In other words the action of $M$ on $V$ is such that 
any element $v \in V$ is contained in a finite dimensional
$M$-stable polynomial submodule  $V_v \subset V$. 

\begin{defi}
Let $M$ be a linear algebraic monoid. A {\em polynomial character}\/  of
$M$ is a  multiplicative morphism  $\rho: M\rightarrow \Bbbk$, such that
$\rho(1)=1$. We will denote as $\mathcal X(M)$ the abstract monoid of all
characters of $M$.

We say that the character $\rho$ is trivial if it only
takes the value $1 \in \Bbbk$.
\end{defi}

\begin{example}
\noindent (1) In the case that $M$ is the monoid of 
$n\times n$ matrices with coefficients in
$\Bbbk$, $\det :M\rightarrow \Bbbk$ is a character. 

\noindent (2) In the case of affine algebraic $\Bbbk$-groups, 
the above definition of character coincides with the usual one. 
Moreover, if $\rho$ is a character of $M$ and $a \in G=G(M)$, then
$\rho(a) \neq 0$: otherwise
$1=\rho(1)=\rho(aa^{-1})=\rho(a)\rho(a^{-1})=0$.
Then, the restriction of a character of $M$ yields 
a character of $G$.
\end{example}

\begin{obse}
\label{obse:varias}
\noindent (1)  The description of $\operatorname {GL}_n$ as the set of
matrices with non zero determinant can be generalized for general
affine monoids. Indeed, if $M$ is an arbitrary affine monoid, using
Theorem \ref{theo:closedinm} we can suppose that $M$ 
is a closed submonoid of ${\operatorname M}_n(\Bbbk)$ for some $n$. Hence the
restriction of $\operatorname{det}$ to $M$ is a character, and the
equality $G(M)= \operatorname {GL}_n(\Bbbk) \cap M$ reads as 
\[
G(M)=M_{\operatorname{det}}= \bigl\{m \in M: \operatorname{det}(m)
\neq 0\bigr\}\,. 
\]

\noindent (2) Observe that if $0 \in M$ and $\rho$ is a non-trivial
polynomial character, then $\rho(0)=\rho(0m)=\rho(0)\rho(m)$ for all
$m\in M$,  and thus $\rho(0)=0$.

\noindent (3) It is clear that the polynomial characters of $M$
correspond bijectively with the polynomial module structures on $\Bbbk$.
If $\lambda$ is a polynomial character of $M$, 
it is easy to show that for a
polynomial representation $V$ of $M$ an element $v \in V$ is a
semi-invariant of weight $\lambda$  if and only if $\chi_V(v)=v
\otimes \lambda$. 
\end{obse}

In the same manner than for algebraic groups, it is clear that
the category of polynomial $M$-modules is an
abelian tensor category, with $\Bbbk$ as unit for the tensor
product. Moreover, the finite dimensional polynomial $M$-modules, admit
duals in the usal manner. Observe that if $V$ is a
left polynomial $M$-module, then the dual $V^*$ is a right polynomial
$M$-module.

If $V$ is a polynomial $M$-module, then the symmetric algebra built
on $V$, that we denote as $\operatorname {S}(V) = \bigoplus_n
\operatorname {S}^n(V)$, is also a polynomial module.
 
It is important to notice that if $M$ is an affine
algebraic monoid, $\Bbbk [M]$ can be naturally endowed with two
structures of a polynomial $M$-module, one from the left and the other
from the right. In explicit terms, if $f \in \Bbbk[M]$ and we
denote as usual $\Delta (f) =\sum f_1 \otimes f_2$, then if $m \in M$  
we have that $m \cdot f = \sum f_1 f_2(m)$ and $f \cdot m = \sum
f_1(m) f_2$.

In particular,  in a similar manner than for the
case of a group and a subgroup, one can define the induction
functor from the representations of a submonoid to the
representations of the whole monoid, see
\cite{kn:doty_representations} and \cite{kn:rennersurvey} for the
definitions in the case of 
monoids and \cite{kn:nosotros} for the case of groups.

\begin{defi} 
If $N\subset M$ is a submonoid, then the
restriction functor from the polynomial $M$-modules to the polynomial
$N$-modules
admits a right adjoint, called the {\em induction functor}\/ from $N$ to
$M$, and that is denoted as  $\operatorname{Ind}_N^M$. 
\end{defi}

Given a finite dimensional polynomial $N$-module $V$,
$\operatorname{Ind}_N^M(V)$ is obtained 
in the following way  (see \cite{kn:doty_representations}).  Call
$\operatorname{Mor}(M,V)$ the
vector space of all the morphisms of $\Bbbk$-varieties 
and take the subspace $\operatorname{Ind}_N^M(V) \subset 
\operatorname{Mor}(M,V)$ given by all $f \in \operatorname {Mor}(M,V)$
such that $f(nm)=n\cdot f(m)$ for all $n\in N$, $m\in M$. We endow this
subspace with the
polynomial $M$-module structure defined as follows: 
$(m'\cdot f)(m)=f(mm')$,  $m,m'\in M$, $f\in
\operatorname{Ind}_N^M(V)$.

The map $\epsilon_V : \operatorname{Mor}(M,V)\rightarrow  V$, \,
$\epsilon_V(f) = f(1)$ restricts to $\operatorname{Ind}_N^M$, and it is a
morphism  of $N$-modules, that is the counit of the adjunction.

In the case that $V$ is infinite dimensional, we proceed in the same
manner than for affine algebraic groups, first we have to take inside of
$\operatorname{Mor}(M,V)$ the subspace 
$\operatorname{Mor}_{\text {fin}}(M,V)=\bigl\{f: M \rightarrow V
\mathrel{:} \dim \langle f(M) \rangle_\Bbbk< \infty \bigr\}$ and
then take in  $\operatorname{Mor}_{\text {fin}}(M,V)$ the
$N$-equivariant maps as before.

Similarly than for algebraic groups, 
one can also define 
\[
\operatorname{Ind}_N^M(V) = \bigl\{ \sum f_i \otimes v_i ,\,
f_i \in \Bbbk[M], v_i \in V : \, \forall \, n \in N,
\sum f_i \cdot n \otimes v_i = \sum f_i \otimes n \cdot v_i \bigr\}
\]
 with the 
structure of polynomial $M$-module 
given by the action of $M$ on $\Bbbk[M]$
{\em on the left}.
 
For future reference we write down the universal property of the
induction functor. 
Given a polynomial $N$-module $V$, a polynomial $M$-module $U$ and a
homorphism of polynomial  $N$-modules  $\varphi: U\rightarrow V$,
then there exists an unique homomorphism
$\widetilde{\varphi}:U\rightarrow  \operatorname{Ind}_N^M(V)$ such
that $\varphi=\epsilon_V\circ \widetilde{\varphi}$.

It is easy to prove the usual properties of the induction functor, for
example the transitivity. Given a chain of closed submonoids 
$S\subset N\subset M$ and a
polynomial $S$-module $V$, then $\operatorname{Ind}_S^M(V) \cong
\operatorname{Ind}_N^M\bigl(\operatorname{Ind}_S^N(V)\bigr)$.    

Moreover, the induction functor is left exact, since it is a right
adjoint of the restriction functor.

The following considerations will be used later.

Assume that $V$ is a polynomial $M$-module and call $G= G(M)$ the group of
invertible elements. It is clear that by restricting the action from
$M$ to $G$ we can view $V$ as a rational $G$-module. If $\chi$ is
the associated $\Bbbk[M]$-comodule structure on $V$, then the
associated $\Bbbk[G]$-comodule structure, that we denote as
$\widehat{\chi}$, is defined by the commutativity of the diagram below:

\[
\xymatrix{ 
V\ar[r]^-{\chi}\ar[dr]_{\widehat{\chi}}&V \otimes \Bbbk[M]
\ar@{^{(}->}[d]^{id \otimes \pi}\\
 & V \otimes \Bbbk[G]
}
\]

where $\pi: \Bbbk[M] \rightarrow \Bbbk[G]$ is the restriction map that
can be viewed as an inclusion.

The following definition is relevant to pin down the difference between
representations of $G$ and of $M$.

\begin{defi} In the situation above a rational representation $\theta$
  of $G$
  is said to be {\em polynomial} (with respect to $M$), if there
  exists a polynomial representation $\chi$ of $M$ with the property
  that $\widehat{\chi}= \theta$. 
\end{defi}

\begin{obse} 

\noindent (1) If $\theta: V \rightarrow V \otimes
  \Bbbk[G]$ is a polynomial representation, the polyonomial
  representation $\chi$ of $M$ that yields $\theta$ is unique, since
  $(\operatorname {id} \otimes \pi) \chi = \theta$ and $\pi$ is
  injective.

\noindent (2) An action $G \times V \rightarrow V$ is 
 a polynomial  representation of $G$ if and only if it 
can be extended to a polynomial
 action of $M$, in other words, if and only if one can find an action
$M \times V
 \rightarrow V$ such that the diagram below commutes.
\[
\xymatrix{ 
G \times V\ar[r]\ar@{^{(}->}[d] & V \\
M \times V \ar[ur] & 
}
\]

Indeed, in the notations above, if $\theta(v)=\sum v_0 \otimes v_1|_{_G}$
with $v_1 \in \Bbbk[M]$, then the action $g\cdot v = \sum v_0 v_1(g)$
can be extended to the polynomial action $m \cdot v= \sum v_0 v_1(m)$.

\noindent (3) Equivalently, if $V$ is finite dimensional the morphism
$G \rightarrow \operatorname{GL}(V)$ produces a polynomial
representation of $G$ if and only if there is a multiplicative morphism 
$M \rightarrow \operatorname {End}(V)$ that restricts to the original
action of $G$. 

In other words if the diagram below is commutative:

\[
\xymatrix{ 
G\ar[r]\ar@{^{(}->}[d] & \operatorname{GL}(V)\ar@{^{(}->}[d]\\
M \ar[r] & \operatorname{End}(V)
}
\]

\noindent (4) If $V$ is a polynomial $G$-module, then it is clear that
$S(V)$ is also a polynomial $G$-module. 
\end{obse}

In the case that $G$ is an affine algebraic group and $M$ is a monoid
equipped with a polynomial 
character $\rho:M \rightarrow \Bbbk$ with the property
that $M_{\rho} = G=G(M)$ (see Observation \ref{obse:varias}), more can be said concerning the
relationship between the rational and polynomial representations of $G$.
\begin{lema} \label{lema:tonto} In the situation above, if $\chi:V \rightarrow V \otimes
  \Bbbk[G]$ is a rational finite dimensional representation of $G$,
  there is an exponent $n \geq 0$ with the property that the rational
  representation of $G$, $\rho^n|_{_G}\,\chi: V \rightarrow V \otimes
  \Bbbk[G]$ is polynomial.
\end{lema}
\dem The result follows immediately from the fact that
$\Bbbk[M]_\rho=\Bbbk[G]$.
\qqed

\begin{obse} 
\noindent (1) Notice that in Lemma \ref{lema:tonto} we have used the
  $\Bbbk[G]$-module structure of $V \otimes \Bbbk[G]$ given by
  multiplication on the second tensorand. 

\noindent (2) Recall that if
  $- \cdot -: G \times V \rightarrow V$ is a 
  rational representation of $G$ and $\lambda$ is a character of $G$, 
we can define a new
  rational finite dimensional representation of $G$ that we denote as
$- \cdot_\lambda -: G \times V \rightarrow V$ by the formula $a
  \cdot_\lambda v= \lambda(a) a \cdot v$ for all $a \in G$ and $v \in V$. 
Using this notation, another
  formulation of the above lemma is the following: if 
$- \cdot -: G \times V \rightarrow V$ is a finite dimensional
  rational representation of $G$, then for some positive integer $n$
  the representation $- \cdot_{{\rho^n|_{_G}}}-: G \times V \rightarrow V$
  is polynomial.    

\noindent (3) In the situation above, 
assume that $n_V$ is the minimal exponent with the
  property that after twisting the orginal action of $G$ on $V$ by
  $\rho^{n_V}$ we obtain a polynomial action -- notice that the
  exponent $n_V$ may be negative. Then, $n_{V
  \otimes W}=n_V+n_W$ and $n_{V \oplus W}=\operatorname {max}\{n_V,n_W\}$.

\end{obse}
Assume that we have an element $\lambda \in \mathcal X(M)$ and $v \in
V$ such that $\lambda \chi_{_M}(v)= v \otimes 1$. Then $v$ is a
semi-invariant for the action of $G$ on $V$ of weight $\lambda|_{_G}^{-1}$.
Indeed, if $\chi_M(v)=\sum v_0 \otimes v_1$, then $\sum v_0 \otimes
\lambda v_1= v \otimes 1$. Evaluating at $g \in G$ we have that
$\lambda(g) \sum v_1(g)v_0= v$ or equivalently that 
$\lambda(g)(g \cdot v)= v$.

\section{Representations and characters of  reductive monoids}

The representation theory of reductive groups is well understood, and
its principal features can be 
extended to the representations of reductive monoids. This extension
was first performed  by 
S. Doty in \cite{kn:doty_representations} for a special class of reductive
monoids, namely the ones with unidimensional center, and then by
Renner in full generality (see \cite{kn:rennersurvey}). 

In this section we present the basic results about the
representation theory of reductive monoids, that will be needed in
the rest of the paper. We will use (as in \cite{kn:am}) the methods of
the classification theory  of reductive monoids. In some places in
order to make the presentation clearer and consistent with our
exposition, some of the proofs of the 
needed results will be presented in some detail. the results we present
in this section are valid in arbitrary characteristic.

\begin{obse}
In accordance with the results of \cite{kn:am}, if $G$ is a semisimple
group the unique algebraic monoid $M$ with unit group  
$G(M)=G$ is $M=G$. Hence in what follows we assume that $G$
is a reductive group with connected center of dimension greater than
$0$.    
\end{obse}

\begin{theo}[Extension principle, see \cite{kn:rennersurvey}]
\label{theo:extension} 
Let $M$ be
a reductive monoid with  unit group $G$\,,  and let $T$ be a maximal torus of $G$ with closure
$\overline{T}$ in $M$. If $S$ is an algebraic
monoid and $\alpha:  G\rightarrow S$, 
$\beta: \overline{T}\rightarrow S$ are morphisms of algebraic monoids such that
$\alpha|_T=\beta|_T$. Then there exists a morphism of monoids
$\varphi:M\rightarrow S$ such that $\varphi|_{_G}=\alpha$ and
$\varphi|_{\overline{T}}=\beta$.
\end{theo}

\dem
Let $X_1,\ldots , X_r$ be the $(G\times G)$-stable
irreducible divisors of $M$ and call $\mathcal W$ 
the union of the corresponding $(B\times 
B^-)$-open orbits (see \cite{kn:knop}). Let
$\mathcal U$ be the    
union of the $(B\times B^-)$-orbits of codimension smaller or equal than
one. As  the complement of $\mathcal
U$ in $M$ is the union of the $(B\times B^-)$-orbits of codimesion
greater or equal than $2$,  $\mathcal U$ is open and 
every morphism of varieties $\psi:\mathcal U\rightarrow S$ extends to
a  morphism $\varphi : M\rightarrow S$. By dimensional reasons, it
follows that  the $(B\times
B^-)$-stable subset $\mathcal A=BB^-\cup
\mathcal W\subset \mathcal U$  is an open
subset. 

Let $e_i\in X_i\cap \overline{T}$ be an idempotent element such that
$(B\times B^-)\cdot e_i$ is the open 
$(B\times B^-)$-orbit of $X_i$ (see \cite{kn:am}). As $e_i\in
\overline{T}\setminus T$, it 
follows that if we consider the action $t\cdot x=tx$, $t\in T$, $x\in
M$, then the isotropy group $T_{e_i}$ is  of dimension greater or equal than
one.  On the
other hand, $B\times B^-=TU\times TU^-$, and thus 
$(UT\times U^-)\cdot e_i=(B\times B^-)\cdot e_i$ is open in
$X_i$. As $\codim X_i=1$,  
\[
1=\dim (UT\times U^-)_{e_i}=\dim (U\times U^-)_{e_i}+\dim T_{e_i}\geq 1\,. 
\]

In follows that  $\dim(U\times U^-)_{e_i}=0$, and that if $v\in X_i\cap
\mathcal W$, then 
$\dim (U\times U^-)_{v}=0$. Thus, the morphism 
$
\varphi:(U\times U^-)\times
(\mathcal A \cap \overline{T})\rightarrow \mathcal A=(U\times
U^-)\cdot (\mathcal A\cap \overline{T})$,
  given by $\bigl((a,b)),v\bigr)\mapsto avb$, where
  $\bigl((a,b)),v\bigr)\in (U\times U^-)\times
(\mathcal A \cap \overline{T})$, is a
birational morphism with finite fibres. As both $(U\times U^-)\times
(\mathcal  A \cap \overline{T})$ and $\mathcal A$ are irreducible and
 $\mathcal A$ is normal (see \cite{kn:am}), it follows from Zariski
 main theorem (\cite[Ch III, \S 11]{kn:hartshorne}) that
$\varphi$ is an isomorphism. In particular, $(U\times U^-)_v=\{1\}$ for all
$v\in \mathcal A$.

Define $\gamma:\mathcal A\rightarrow S$ by
$\gamma(m)=\alpha(a)\beta(v)\alpha(b)$, if $m=avb$, $a\in U$,
$b\in U^-$, 
$v\in \mathcal A\cap \overline{T}$. 
Then $\gamma|_{\mathcal A \cap G}=\alpha |_{\mathcal A \cap G} $, and
thus we have a morphism $\psi:\mathcal U\cup G \rightarrow S$
 such that
$\psi(mm')=\psi(m)\psi(m')$, for all $m,m'\in \mathcal W$. It
follows by continuity that the extension of $\psi$ is a morphism of
monoids that verifies the required properties.
\qed

\begin{defi}
Let $M$ be a reductive
monoid of unit group $G$ and call $\overline{T}$ the closure of
the maximal torus $T$. We define the set of {\em polynomial 
weights}\/ of $M$, denoted as $\mathcal X(\overline{T})$, as the set of 
all multiplicative morphisms of algebraic monoids from $\overline{T}$ 
to $\Bbbk$. Notice that  $\mathcal X(\overline{T})$ can be interpreted
as the submonoid of $\mathcal X(T)$ consisting of the 
characters of $T$ that extend to $\overline{T}$. We define  the {\em
dominant polynomial weights}\/ as the intersection
$\mathcal X_+(\overline{T})=\mathcal X_+(T)\cap \mathcal X
(\overline{T})$. 
\end{defi}

\begin{theo}
Let $M$ be a reductive monoid, and let $V$ be a rational finite
dimensional $G$-module 
all whose $T$-weights are polynomial, then $V$ is a polynomial $G$-module.
\end{theo}

\dem
Let $\{v_i\}_{i\in I}$ be a basis of $V$ consisting of $T$-weight
vectors. If we write $g\cdot v_i=\sum c_{ij}(g)v_j$, with $c_{ij} \in 
\Bbbk[G]$, then the functions $c_{ij}|_{T}$ extend to $\overline{T}$
for all $i,j\in I$. Then, it follows from the  extension principle
(Theorem \ref{theo:extension}) that $c_{ij}$ extend to $M$ for all
$i,j\in I$. 
\qqed

Consider $\lambda\in \mathcal X(\overline{T})$, and consider the one
dimensional $T$-module $\Bbbk_\lambda$, where as usual $\Bbbk_\lambda$
is the base field with the action $t\cdot 
\alpha=\lambda(t)\alpha$, $t\in T$, $\alpha\in \Bbbk$. Then
$\Bbbk_\lambda$ is a $\overline{T}$-mod, and if $\overline{B}$ is the
closure of $B$ in $M$, this $\overline{T}$-module structure extends to a 
structure of $\overline{B}$-module.

\begin{theo}
Let $M$ be a  reductive monoid and consider $\lambda\in
\mathcal X(\overline{T})$. Then
$\operatorname{Ind}^M_{\overline{B}}\Bbbk_\lambda$ is non-zero if and only
if $\lambda\in \mathcal X_+(\overline{T})$.

In particular, if
$L(\lambda)=\operatorname{soc}_G\operatorname{Ind}^{^G}_{\overline B}
\Bbbk_\lambda $, then the set of isomorphism classes of simple
polynomial 
$M$-modules is $\bigr\{L(\lambda): \lambda\in \mathcal
X_+(\overline{T})\bigr\}$. 
\qedennonce
\end{theo}

Next we study the properties of the abstract monoid $\mathcal
X_+(\overline{T})$.

\begin{theo}
\label{thm:xi(t)}
In the above situation
\[
\mathcal X_+(\overline{T})=
\left\{ \lambda\in \mathcal X_+(T)\ |\ \exists f\in
\Bbbk[M]^{(B\times B^-)}\,\right\}=\mathcal C^\vee\cap \mathcal X_+(T). 
\] 

As usual, $\Bbbk[M]^{(B\times B^-)}\,=\left\{f \in \Bbbk[M]:\ (a,b)\cdot
f=\lambda(a)\lambda^{-1}(b)f\ \forall a,b \in T\right\}\,$.
\end{theo}

\dem 
If $\lambda\in \mathcal X_+(T)$, then there exists a nonzero $f\in
\Bbbk[G]$ such that 
$(a,b)\cdot f=\lambda(a)\lambda^{-1}(b)f$, $(a,b)\in B\times
B^-$ and then, $f(x)=\bigl((x,1)\cdot f\bigr)(1)=\lambda(x)f(1)$
for all $x\in T$. 

Assume also that $f\in \Bbbk[M]$. If
$f(1)=0$, then for $x=ab\in BB^-$ we have $f(x)=\bigl((a,b^{-1}\bigr)\cdot
f)(1)=\lambda(a)\lambda^{-1}(b^{-1})f(1)=0$, and it follows that
$f=0$. Then $f(1)\neq 0$, and in this situation we can extend 
$\lambda$ to all $\overline{T}$ using the formula
$\lambda(x)=\frac{f(x)}{f(1)}$, and thus $\lambda \in
\mathcal X(\overline{T})$. Hence,  $\lambda\in \mathcal X_+(\overline{T})$. 

On the other hand, let $\lambda\in \mathcal X_+(\overline{T})$, and
consider $f\in 
\Bbbk[G]$ such that $(a,b)\cdot f=\lambda(a)\lambda^{-1}(b)f$ for all
$(a,b)\in B\times B^-$; we have
$f(z)=f(1z)=\lambda(z)f(1)$ for $z\in T$.
Let $g:\overline{T}\rightarrow \Bbbk$ be the morphism given by $
g(z)=\lambda(z)f(1)$. As $g(t)=f(t)$ for all
$t\in T$,  by the extension principle (Theorem \ref{theo:extension})
there 
exists  $h\in \Bbbk[M]$ such 
that $h|_{\overline{T}}=g$ and $h|_{_G}=f$. As
$\bigl((a,b)\cdot h\bigr)(x)=\lambda(a)\lambda^{-1}(b)h(x)$ for all
$x\in G\subset 
M$ and $(a,b)\in B\times B^-$, it follows by continuity that $h$ is a
weight vector of 
weight $(\lambda,-\lambda)$ for $B\times B^-$.\qed

\begin{theo}
\label{thm:overline(t)}
In the situation above the monoid $\mathcal X_+(\overline{T}) $ is
saturated, i.e., $n\lambda\in  \mathcal X_+(\overline{T})$, if and only if
$\lambda\in \mathcal X_+(\overline{T})$. Moreover, if $M$ has a zero, then
$\mathcal X_+(\overline{T})$ is an ideal of $\mathcal X_+(T)$, i.\
e.\ if $\mu\in 
\mathcal X_+(\overline{T})$ and  $\lambda\in \mathcal X_+(T)$ is such
that $\lambda\leq 
\mu$, then $\lambda\in \mathcal X_+(\overline{T})$.
\qedennonce
\end{theo}

We finish this section with an existence result that will allow us to
apply Cayley--Hilbert's methods  in order to  construct a finite
number of generators for the algebra of invariants of a linear
representation. 

\begin{lema}
Let $G$ be a reductive non semisimple group and $\lambda\in \mathcal X(G)$ a
character. Then there exists an algebraic monoid with zero
$M$ that has $G$ as a unit group and with the property that $G=M_\lambda$.
\end{lema}

\dem
Let $w=\sum_{i=1}^l \omega_i$ and $\mathcal C$ the cone generated by
$(-w,\lambda)$ and all the coroots. It is clear that $\bigl(\mathcal C,
\{\alpha_1^\vee,\dots, \alpha_l^\vee\}\bigr)$ is a colored cone
associated to an affine variety $M$. Indeed, it is a striclty
convex rational cone, and 
$(0,\lambda)=(-w,\lambda)+(\frac{1}{2}\sum \alpha_i,0)$; hence,
it follows from
the classification theory of reductive monoids that $M$ is  a reductive
monoid with zero and unit group $G$ (see \cite{kn:am}).

We deduce from the construction that $\lambda$ is a polynomial
character for $M$, since $\lambda\in \mathcal C^\vee\cap \mathcal
X(G)$. Hence $M_\lambda$ is a $(G\times G)$-stable divisor, 
and thus it is the unique $(G\times G)$-stable divisor of $M$,
 corresponding to the edge of $\mathcal C$ generated by $(-w,\lambda)$.
\qed

\section{Generalized Cayley's $\Omega$-processes}
\label{section:process}

Let $M$ be an affine algebraic monoid and assume that $\lambda$ is a
non-trivial character of $M$.

\begin{defi}
An {\em $\Omega$-process} (associated to $\lambda$) is a non zero
linear operator
  $\Omega:\Bbbk[M]\to \Bbbk[M]$ such that:

\noindent $\Omega(f\cdot
  m)=\lambda(m)\Omega(f)\cdot m$ ; $\Omega(m\cdot
  f)=\lambda(m)m\cdot \Omega(f)$ for all $f\in \Bbbk[M]$ and $m\in M$.
\end{defi}

We will usually abbreviate and call an $\Omega$-process simply a process.
 
In classical nomenclature the above definition is called the ``first
rule of a process''.
 
\begin{obse}
 From the above definition one easily concludes that if
  $r$ is a positive integer, then $\Omega^r(f\cdot
  m)=\lambda^r(m)\Omega(f)\cdot m$ ; $\Omega^r(m\cdot
  f)=\lambda^r(m)m\cdot \Omega(f)$ for all $f\in \Bbbk[M]$ and $m\in
  M$. In other words, if $\Omega$ is a process associated to  
$\lambda$, then $\Omega^r$ is a process associated to
  $\lambda^r$.
 \end{obse}

\begin{example}[The classical $\Omega$-process]
\label{exam:clasico}
Let $M=\operatorname{M}_n(\Bbbk)$. Then the diferential operator of
order $n$, $\Omega: \Bbbk[x_{11},\dots, x_{nn}]\to \Bbbk[x_{11},\dots,
x_{nn}]$, given by
\begin{equation}
\label{eqn:clasico}
\Omega(f)=\sum_{\sigma\in \mathcal
  S_n}\operatorname{sg}(\sigma)\frac{\partial^n}{\partial
  x_{1\sigma(1)}\cdots \partial x_{n\sigma(n)}}(f)
\end{equation}
is an $\Omega$-process. See \cite{kn:Sturm} for a modern presentation.
\end{example}

\begin{lema}
With the notations above, let $\mu\in \mathcal X(M)$ be a polynomial
character of the monoid, and $g\in \Bbbk[M]$ a right
    $M$-semi-invariant of weight $\mu$. Then  for all $m \in M$
\begin{equation}
\label{eqn:basic}
\Omega \bigl((f\cdot m )g\bigr)\mu(m)= \Omega\bigl((fg)\cdot m\bigr)= 
\lambda(m)\Omega(fg)\cdot m \,,
\end{equation}

Similarly, if $h$ is a left $M$-semi-invariant of weight $\mu$, then for
all $m\in M$

\begin{equation}
\label{eqn:basic2}
\Omega\bigl((m \cdot f)h\bigr)\mu(m)= \Omega\bigl(m \cdot (fh)\bigr)= 
\lambda(m)m \cdot \Omega(fh)\,;
\end{equation}

In particular,

\begin{equation}
\label{eqn:basic3}
\Omega(f\lambda^s)\cdot m\,= 
\Omega\bigl((f\cdot m) \lambda^s\bigr)\lambda^{s-1}(m)\quad \quad \forall s\geq 1\,, 
\end{equation}

and

\begin{equation}
\label{eqn:basic4}
m\cdot \Omega(f\lambda^s)= 
\Omega\bigl((m \cdot f) \lambda^s\bigr)\lambda^{s-1}(m)
\quad \quad \forall s\geq 1\,, 
\end{equation}
\end{lema}

\dem
It follows easily from the definition of a process.
\qed

\begin{lema}
In the notations above,
\[\Omega(\lambda^s)\cdot m\, = 
\Omega(\lambda^s)\lambda^{s-1}(m)
\]
\[ m \cdot \Omega(\lambda^s) = 
\Omega(\lambda^s)\lambda^{s-1}(m)  
\] 
 
Moreover, $\Omega(\lambda^s)=\Omega(\lambda^s)(1)\lambda^{s-1}=
\alpha_{s}\lambda^{s-1}$ with $\alpha_{s}= \Omega(\lambda^s)(1) \in \Bbbk$.

\end{lema}

\dem
 Just make the substitution   $f=1$ in equations
 (\ref{eqn:basic}) and (\ref{eqn:basic2}).

As for the second assertion, evaluate the above equality at $1\in M$. 
\qed

\begin{lema}
Let $\mu\in \mathcal X(M)$ be a polynomial character,  and
$g\in\Bbbk[M]$ ($h\in\Bbbk[M]$\,)
 a right (left) $M$-semi-invariant of weight $\mu$, then 
\[
\Omega^r\bigl((f\cdot m )g\bigr)\mu^s(m)= \Omega^r\bigl((fg )\cdot m\bigr)= 
\lambda^r(m)\Omega^r(fg)\cdot m 
\]
and
\[
\Omega^r\bigl((m \cdot f)h\bigr)\mu^s(m)= \Omega^r\bigl(m \cdot (f\mu^s )\bigr)= 
\lambda^r(m)m \cdot \Omega^r(fg).
\]

In particular,
\[
\lambda^r(m)m \cdot \Omega^r(\mu^s) =
\Omega^r(\mu^s)\mu^s(m)=  
\lambda^r(m)\Omega^r(\mu^s)\cdot m \,;
\]
\[
\Omega^r(\mu^s)(1)\mu^s(m)=  
\lambda^r(m)\Omega^r(\mu^s)(m) \,;
\] 

\noindent in other words: $\Omega^r(\mu^s)(1)\mu^s =  
\lambda^r\Omega^r(\mu^s)$.

\end{lema}

\dem 
Recall that $\Omega^r$ is a process associated to $\lambda^r$, and
then apply equations (\ref{eqn:basic})--(\ref{eqn:basic4}).
\qed

\begin{obse} \label{obse:proper}
In particular, if in the above lemma we consider $\mu=\lambda$, we obtain that
$\Omega^r(\lambda^s)(1)\lambda^s= \Omega^r(\lambda^s)\lambda^r$,
and then if we call $\Omega^r(\lambda^s)(1)=\alpha_{r,s}$ -- notice
that in the notations used before $\alpha_s=\alpha_{1,s}$ --, we deduce
that 
\[
\Omega^r(\lambda^s)\lambda^r=\alpha_{r,s}\lambda^s\,. 
\]

From $\Omega(\lambda^s)=\alpha_s \lambda^{s-1}$,
we deduce that $\Omega^2(\lambda^s)=\alpha_s \alpha_{s-1}\lambda^{s-2}$, if $s
\geq 2$. Then $\Omega^r(\lambda^s)=\alpha_s \alpha_{s-1} \ldots
\alpha_{r-s} \lambda^{r-s}$ if $r \geq s$ and if 
we call $c_{s}=\alpha_s \alpha_{s-1}\ldots \alpha_{1}$, then
$\Omega^s(\lambda^s)=c_s \in \Bbbk$.
  \end{obse}

\begin{defi} 
In the situation above we say that the $\Omega$-process is
  {\em proper}\/ if the elements $\alpha_s \in \Bbbk$ are not zero
  for all $s=1,2,\dots $. Recall that 
$\Omega(\lambda) =\Omega(\lambda)(1)=\alpha_1, 
\Omega(\lambda^2)(1)=\alpha_2, \dots ,
\Omega(\lambda^s)(1)=\alpha_s$ for all $s \geq 1$.
\end{defi}

\begin{obse}
If we write $f\cdot m=\sum f_1(m)f_2$, we have 
\[
\Omega^r\bigl((f\cdot m)g\bigr)\mu^s(m)=
\bigl(\sum f_1(m)\Omega^r
(f_2g)\bigr)\mu^{s}(m)=\lambda^r(m)\Omega^r(fg)\cdot m\,. 
\]

Evaluating at $0$ we have that
\[\lambda^r(m)
\Omega^r(fg)(0)=\bigl(\sum f_1(m)\Omega^r
(f_2 g)(0)\bigr)\mu^{s}(m)
\]

i.e.,
\[ \Omega^r(fg)(0) \lambda^r = 
\bigl(\sum f_1\Omega^r(f_2g)(0)\bigr)\mu^s\,.
\]

Similarly 
\[
\Omega^r\bigl((m \cdot f)h\bigr)\mu^s(m)=
\bigl(\sum \Omega^r
(f_1h)f_2(m)\bigr)\mu^{s}(m)=\lambda^r(m)m \cdot \Omega^r(fh) \,. 
\]

and

\[ \Omega^r(fh)(0) \lambda^r = 
\bigl(\sum \Omega^r(f_1h)(0)f_2\bigr)\mu^s\,.
\]

In particular if $g = h = \lambda = \mu$ we have that

\begin{equation}
\Omega^r(f\lambda^s)(0)\lambda^r=\lambda^{s}\sum 
\Omega^r(f_1\lambda^s)(0)f_2= \lambda^s
\sum f_1 \Omega^r(f_2\lambda^s)(0)\,.
\end{equation}

Moreover, if $r=s=1$, 
\begin{equation}
\Omega(f\lambda)(0)=\sum f_1
\Omega(f_2\lambda)(0)= \sum \Omega(f_1\lambda)(0)f_2\,.
\end{equation}
\end{obse}

\section{Existence of Cayley's $\Omega$-processes}

In this section we describe in the case of a reductive monoid, all the
$\Omega$-processes associated to a 
polynomial character $\lambda\in \mathcal X(M)$. 
In order to deduce the general formula for $\Omega$, first we apply
the forumul\ae\ above, in order to calculate the value of the process
on semi-invariants. 

\begin{theo}
\label{theo:semininv}
Let $\mu\in \mathcal X(\overline{T})$ and $f\in\Bbbk[M]$ be a $(B\times
B^-)$-semi-invariant of weight $(\mu,-\mu)$. Then $\Omega(f)$ is a
$(B\times
B^-)$-semi-invariant of weight $(\mu -\lambda,\lambda -\mu)$. 
\end{theo}

\dem
In the same manner that before, from the definition of a
process we deduce that for all $b_1\in B$ and $b_2\in B^-$,
\begin{equation}
\label{eqn:semi-inv}
\begin{array}{ll}
(b_1,b_2)\cdot \Omega(f)\lambda(b_1)\lambda(b_2^{-1})= &
(b_1\cdot\Omega( f)\cdot b_2^{-1}) \lambda(b_1)\lambda(b_2^{-1})=\\ 
& \Omega\bigl( (b_1,b_2)\cdot f\bigr)= 
 \mu(b_1)\mu(b_2^{-1})\Omega(f)\,.
\end{array}
\end{equation}
\qed

\begin{coro}
\label{coro:cuando0}
Let $\Omega$ be a process for $M$ associated to a polynomial character
$\lambda$, and let $f\in \Bbbk[M]$ be a semi-invariant of weight
$(\mu,-\mu)$. If $\mu-\lambda\notin 
\mathcal C^\vee$, then $\Omega(f)=0$.
\end{coro}

\dem
Since $\operatorname{char}\Bbbk=0$, we can decompose 
$\Bbbk[M]=\bigoplus_{\mu\in \mathcal X_+(\overline{T})}
V_\mu\otimes V^*_\mu$, where $V_{\mu}$ is  the
irreducible $G$-module of maximal weight $\mu$ for $B$ -- recall that
accodingly, $V^*_\mu$ is the irreducible $G$-module of maximal weight
$-\mu$ for $B^-$ (see Theorem \ref{thm:xi(t)} and
\cite{kn:am}). Hence, if $\mu-\lambda\notin 
\mathcal C^\vee$ it follows that there is no non-trivial semi-invariant
of weight $\mu-\lambda$ in $\Bbbk[M]$.
\qed

\begin{theo}\label{theo:semininv2}
Let $M$ be a reductive monoid,  $\lambda\in \mathcal X(M)$ and
$\Omega$ a process associated to $\lambda$. Then
$\Omega=\bigoplus_{\mu\in \mathcal X_+(\overline{T})}
  \frac{a_\mu}{\lambda}\operatorname{Id}_\mu$, where
  $\operatorname{Id}_\mu$ is the identity map of $V_\mu\otimes
  V_\mu^*$, and $a_\mu\in \Bbbk $. Moreover, if $\mu-\lambda \notin
  \mathcal X_+(\overline{T})$, then $a_\mu=0$. 

In particular, $\Omega$ is proper if and only if $a_{\lambda^s}\neq 0$ for
all $s\geq 0$.
\end{theo}

\dem From the definition of a process, it is clear that $\Omega$
stabilizes the $G\times G$-submodules of $M$, and hence $\Omega$  is
characterized by its restrictions to $V_\mu\otimes V_\mu^*$, $\mu\in
\mathcal X_+(\overline{T})$. 

If $f_{\mu}\in V_\mu\otimes V^*_\mu$ is a $(B\times
B^-)$-semi-invariant, then it follows from Theorem
\ref{theo:semininv} that 
$\Omega(f_\mu)$ is also a semi-invariant, of weight
$(\mu-\lambda,\lambda-\mu)$. Hence, $\Omega(f_\mu)=a_\mu
\frac{f_\mu}{\lambda}$, where $a_\mu\in\Bbbk$ is necessarily zero if
$\mu-\lambda\notin \mathcal X_+(\overline{T})$.

Again from the definition of a process, it  follows  that if $s,t\in
G$, then  
\[
\begin{split}
\Omega\bigl((s,t)\cdot f_\mu\bigr)= & \Omega\bigl( s\cdot
f_\mu\cdot t^{-1} \bigr)= \lambda(t^{-1})\lambda(s)
\bigl(s\cdot\Omega(f_\mu)\cdot  t^{-1}\bigr) =\\
& a_\mu\lambda(t^{-1})\lambda(s)\frac{ s \cdot f_\mu\cdot
  t^{-1}}{s \cdot \lambda\cdot t^{-1}}=a_\mu \frac{(s,t)\cdot
  f_\mu}{\lambda}\,.    
\end{split}
\]

It follows from Corollary \ref{coro:cuando0} that if
$\mu-\lambda\notin \mathcal X_+(\overline{T})$, then $a_\mu=0$.

Conversely, it is easy to show that any election of $a_\mu\in \Bbbk$, 
with 
$\mu- \lambda \in \mathcal X_+(\overline{T})$, is adequate in order to
have an 
$\Omega$-process.

  Finally, the assertion about properness of an $\Omega$-process follows
  directly from the above definition.
\qqed

\begin{example}
\label{exam:clasico2}
Consider the classical $\Omega$-process (see Example
\ref{exam:clasico}, equation (\ref{eqn:clasico})). Recall that for
$\operatorname{GL}_n(\Bbbk)$ the
fundamental weights are given by 
\[
\omega_i\left(
\begin{smallmatrix}
a_{11} & &0\\
&\ddots &\\
0 & &a_{nn}
\end{smallmatrix}
\right)=a_{n-i+1,n-i+1}\dots a_{nn}\,,
\]
with respective eigenvectors in $k[x_{11},\dots, x_{nn}]$ given by 
\[
f_i=\operatorname{det} 
\left(
\begin{smallmatrix}
x_{n-i+1,n-1+1} &\cdots &x_{n-1+i,n}\\
\vdots& &\vdots\\
x_{n,n-i+1}&\cdots & x_{nn}
\end{smallmatrix}\,.
\right)
\]

In particular, if
$\omega=\omega^{r_1}_1\dots\omega_{n-1}^{r_{n-1}}\operatorname{det}^r$
is any dominant weight,
then $f=f_1^{r_1}\dots f_{n-1}^{r_{n-1}}\operatorname{det}^{r}$ is an
eigenvector for $\omega$. From Theorem \ref{theo:semininv2} it follows
that
\[
\Omega(f)= a_{\omega}f_1^{r_1}\dots 
f_{n-1}^{r_{n-1}}\operatorname{det}^{r-1}\, ,
\]

and $a_{\omega}= \Omega(f)(1)$.

The computation of $a_{\omega}$ for a general weight $\omega$ seems to
be a non trivial combinatorial problem.  For our purposes, all we need
is to prove that $\alpha_s = \Omega(\operatorname {det}^s))(1) \neq
0$. This will follow from the equality: $\Omega^s(\operatorname
{det}^s)=\alpha_s \alpha_{s-1}\cdots \alpha_1$ -- see
Observation \ref{obse:proper} --, and the fact that $\Omega^s(\operatorname
{det}^s) \neq 0$, see \cite{kn:Sturm}. 

\end{example}

\section{Finite generation of invariants}

In this section we show how to construct from an $\Omega$-process, the standard
toolkit of an invariant theorist: integrals, Reynolds operators, etc.

We show that the existence of a proper process 
guarantees the finite generation of the rings of invariants
corresponding to linear actions. Our proof is basically the same than
the original one due to Hilbert, that appeared in \cite{kn:Hil}. The
required finite generation is guaranteed by the fact that one can use
the process in order to produce quite a large number of
semi-invariants (see Theorem \ref{theo:2rule}).  

The results of the theorem that follows were known in the classical
literature as the ``second rule of a process''.

\begin{theo}
\label{theo:2rule}
Let $V$ be a polynomial $M$-module, and  consider the map $I_{r,s}: V
\rightarrow V$, 
$I_{r,s}(v)=\sum v_0\Omega^r(\lambda^sv_1)(0)$. Then: 

\noindent (1) The map $I_{r,s}$ is a morphism of $M$-modules,
with respect to the original action on $V$ twisted by $\lambda^s$.

\noindent (2) For all $v \in V$,
$\lambda^s\chi\bigl(I_{r,s}(v)\bigr)=\lambda^r\bigl(I_{r,s}(v) \otimes
1\bigr)$. 

\noindent (3) Assume that $V$ is a polynomial $G$-module with
structure $\theta$ and consider it as a polynomial $M$-module with
structure $\chi:V\rightarrow V \otimes \Bbbk[M]$ with
$\theta=(\operatorname {id} \otimes \pi)\chi$. 
If $I_{r,s}$ is as above, then for all $v \in V$, $I_{r,s}(v)$ is a
semi-invariant for the action of $G$ on $V$, of weight 
$\lambda^{r-s}$.
\end{theo}

\dem

\noindent (1) The assertion concerning the fact that $I_{r,s}$ is a morphism of
comodules is equivalent to the commutativity of the diagram below.

\[
\xymatrix@C=1.6cm{ 
V\ar[r]^{I_{r,s}}\ar[d]_{\lambda^{s}\chi}&
V\ar[d]^{\lambda^{s}\chi}\\
V \otimes \Bbbk[M] \ar[r]^{I_{r,s} \otimes \operatorname {id}}& 
V \otimes \Bbbk[M]
}
\]

Then, $\lambda^{s}\chi \bigl(I_{r,s}(v)\bigr) =
\lambda^s\sum v_0\otimes v_1\Omega^r(\lambda^sv_2)(0)=\sum
v_0\Omega(\lambda^sv_1)(0)\otimes \lambda^r$.

Moreover, $(I_{r,s} \otimes \operatorname {id})\bigl(\lambda^s
\chi(v)\bigr)=\sum 
I_{r,s}(v_0) \otimes \lambda^s v_1=\sum v_0 \otimes
\lambda^s\Omega^r(\lambda^sv_1)(0)v_2=\sum v_0 \otimes
\Omega^r(\lambda^sv_1)(0) \lambda^r=\sum v_0\Omega^r(\lambda^sv_1)(0)
\otimes \lambda^r$.

Then, the commutativity of the above diagram follows. 

\noindent (2) The assertion of this part is the first 
of the formul\ae\/ proved above.

\noindent (3) Applying $\operatorname {id} \otimes \pi$ to the
equality just proved we obtain 
that $\lambda^s|_{_G}\,\theta\bigl(I_{r,s}(v)\bigr)= I_{r,s}(v) \otimes
\lambda^r|_{_G}$, and hence we conclude 
that $I_{r,s}$ is a $G$-semi-invariant of weight $\lambda^{r-s}$.
\qed

The definition of integrals for monoids is the same than for groups.

\begin{defi} 
Assume that $M$ is an affine algebraic monoid, a map
  $J:\Bbbk[M] \rightarrow \Bbbk$ is a {\em left normalized integral}\/
  if it 
  is a morphism of $M$-modules with $J(1)=1$ and 
satisfies: $J(f)1=\sum f_1 J(f_2)$. The equality $J(f)1=\sum
  J(f_1)f_2$ characterizes right integrals.
\end{defi}

See \cite{kn:nosotros} for the definition and properties of integrals
and Reynolds operators for algebraic groups. In particular, in
\cite{kn:nosotros} it is shown how to construct Reynolds operators
from integrals. Here we are considering the obvious adaptation of these 
concepts to the context of affine algebraic monoids.

Notice also that if we evaluate at $x$ the equality $J(f)1=\sum f_1
J(f_2)$, we obtain: $J(f)=\sum f_1(x) J(f_2)= \sum J\bigl(f_2 f_1(x)\bigr)=
J(f\cdot x)$.

\begin{theo}
 Let $M$ be an affine algebraic monoid with $0$ and let
  $\lambda$ be a polynomial character of $M$; assume
  that $\Omega$ is a proper process with respect to
$\lambda$. Then, the map $J: \Bbbk[M] \rightarrow
  \Bbbk$, $J(f)=\frac{1}{\Omega(\lambda)}\Omega(f\lambda)(0)$ is a
  two sided normalized integral.

\end{theo}
\dem
The lemma just proved guarantees that $J$ satisfies the above
equality, and clearly $J(1)=1$. The $M$-equivariance of $J$ also
follows from the above.
\qqed

\begin{theo}
With the notations above and assuming that $\Omega$ is a proper
process, ${\mathcal R}_V=\frac{1}{\Omega(\lambda)}I_{1,1}:V\to V$,
is a Reynolds operator for 
the category of polynomial $M$-modules.
\end{theo}
\dem
Indeed, 
$\sum v_0\Omega(\lambda v_1)(0)=I(v)$ is an $M$-invariant.
Now, if $v\in V^M$, then $\chi(v)=v\otimes 1$ and $I_V(v)=
v\Omega(\lambda )(0)/\Omega(\lambda )(0)=v$. 
One can easily show that if $f:V \rightarrow W$ is a morphism of
$M$-modules then the corresponding morphism 
$\widehat{f} = f|_{V^M}: V^M \rightarrow W^M$ 
is compatible with the operators $I_V$ and $I_W$ in the
sense that the following diagram is commutative.

\[
\xymatrix{ 
V\ar[r]^{I_V}\ar[d]_f&V^M\ar[d]^{\widehat{f}}\\
W\ar[r]^{I_W}&W^M 
}
\]
\qqed

The proof of the theorem that follows is only sketched as it is
closely related to the standard methods introduced in \cite{kn:Hil};
see also \cite{kn:Sturm} for a modern presentation in the case of 
$\operatorname {GL}_n$.

\begin{theo}
\label{theo:main}
 Let $M$ be a linear algebraic monoid with $0$ and assume
  that for some polynomial character $\lambda:M \rightarrow \Bbbk$,
  $G=G(M)=M_{\lambda}=\{m \in M: \lambda(m)\neq 0\}$. Assume moreover
  that $M$ admits a proper process with respect to
  $\lambda$. If $V$ is a finite dimensional rational $G$-module, then
  the ring of invariants $\operatorname {S}(V)^G$ is finitely
  generated. 
\end{theo}

\dem Let $I = \langle \operatorname {S}_+(V)^G \rangle \subset
\operatorname {S}(V)$ be the ideal of $\operatorname {S}(V)$ generated
 by the homogeneous invariants of positive degree. Let   $\{\xi_1,\dots
,\xi_t \}  \subset \operatorname {S}_+(V)^G$ be a finite number of
homogeneous ideal
generators of $I$, and call $d_i$ the degrees of the $\xi_i$ for
$i=1,\dots ,t$.

Evidently $\Bbbk[\xi_1,\dots ,\xi_t] \subset \operatorname {S}(V)^G$;
we will prove by induction on $d >0$ that $\operatorname 
{S}_d(V)^G \subset \Bbbk[\xi_1,\dots ,\xi_t]$.

Assume that for all $e < d$,  $\operatorname {S}_e(V)^G \subset
\Bbbk[\xi_1,\dots ,\xi_t]$ and take $\xi \in \operatorname {S}_d(V)^G
\subset I$. Then, we can find homogeneous elements $f_i \in
\operatorname {S}_{e_{i}}(V)$ such that  $\xi = \sum f_i \xi_i$, 
with $d=e_i+d_i$, and $e_i < d$ for $i=1,\dots ,t$.

Call $\theta: V \rightarrow V \otimes \Bbbk[G]$ the comodule structure
on $V$. As we already observed for some $n \geq 0$, $\lambda^n
\theta:V \rightarrow V \otimes \Bbbk[G]$ is polynomial, i.e. it is of
the form $\lambda^n \theta= (\operatorname {id} \otimes \pi)\chi$ for
some $\chi:V \rightarrow V \otimes \Bbbk[M]$, polynomial
$\Bbbk[M]$-comodule structure on $V$.

If we call $\theta_r: \operatorname {S}_r(V) \rightarrow
\operatorname {S}_r(V) \otimes \Bbbk[M]$, the comodule structure
induced by $\theta$ on $\operatorname {S}_r(V)$, then $\lambda^{nr}\theta_r$
is polynomial.

Applying $\theta_d$ to the equality $\xi = \sum f_i \xi_i$ we obtain: 

\[
\xi \otimes 1 =
\sum_{i=1}^t f_{i0} \xi_i \otimes f_{i1}\,. 
\]

Multiplying by $\lambda^{nd}$ we have:  \[\xi \otimes \lambda^{nd} =
\sum_{i=1}^t f_{i0} \xi_i \otimes \lambda^{nd_i}
\lambda^{ne_i}f_{i1}\, ,\]

In other words, for the coaction $\lambda^{ne_i}\theta_{e_i}:
\operatorname {S}_{e_i}(V) \rightarrow \operatorname {S}_{e_i}(V)
\otimes \Bbbk[M]$, that sends
$\lambda^{ne_i}\theta_{e_i}(f_i) = \sum f_{i0} \otimes g_{i1}$ with
$f_{i0} \in \operatorname {S}_{e_i}(V)$ and $g_{i1} \in
\Bbbk[M]$, we have that  
\[\xi \otimes \lambda^{nd} = \sum_{i=1}^t
f_{i0} \xi_i \otimes \lambda^{nd_i} g_{i1}\, .\]

Applying $\Omega^{nd}$, then evaluating at zero and 
recalling that $\Omega^{nd}(\lambda^{nd})=c_{nd} \in \Bbbk^{*}$, we obtain
\[c_{nd}\xi = \sum_{i=1}^t f_{i0} \Omega^{nd}(\lambda^{nd_i}
g_{i1})(0)\xi_i =\sum_{i=1}^t I_{nd,nd_i}(f_i)\xi_i\, ,\] where
$I_{nd,nd_i}$ is as in Theorem \ref{theo:2rule}, where it is proved
that it is 
a semi-invariant of weight $\lambda^{ne_i}$.

As all $\xi, \xi_{i} , i=1,\dots ,t$ are in fact invariants, we
conclude that $ I_{nd,nd_i}(f_i) \in \operatorname {S}_{e_{i}}^G(V)
\subset \Bbbk[\xi_1,\dots,\xi_t]$. Then as $c_{nd} \neq 0$, we have
that $\xi \in
\Bbbk[\xi_1,\dots,\xi_t]$ and the proof is finished.
\qed

\begin{coro}
Let $G$ be a semisimple group and $V$ a finite generated rational
$G$-module. Then $S(V)^G$ is finite generated. 
\end{coro}

\dem
We consider $G'=G\times \Bbbk^*$ that is a reductive non semisimple
group and extend the action of $G$ on $V$ to $G'$  by $(g,z)\cdot
v=g\cdot v$. The result follows directly from Theorem
\ref{theo:main}. 
\qed

\vspace*{1cm}
\noindent{\sc Alvaro Rittatore} \hfill {\sc Walter Ferrer
  Santos}\hspace*{2.5em}\\ 
Facultad de Ciencias \hfill Facultad de Ciencias\hspace*{5em}\\
Universidad de la Rep\'ublica \hfill 
Universidad de la Rep\'ublica\hspace*{1.8em}\\
Igu\'a 4225\hfill Igu\'a 4225\hspace*{9.7em}\\
11400 Montevideo \hfill 11400 Montevideo\hspace*{6em}\\
Uruguay\hfill Uruguay \hspace*{9.8em}\\
e-mail: {\tt alvaro@cmat.edu.uy}\hfill e-mail: {\tt wrferrer@cmat.edu.uy}

\end{document}